\documentclass[12pt,a4paper,twoside]{article}

\newcommand{\pageformat}[6]{\setlength{\hoffset}{-1in}
                  \setlength{\voffset}{-1in}
                  \addtolength{\hoffset}{#5}
                            \addtolength{\voffset}{#6}
                            \setlength{\oddsidemargin}{#1}
                            \setlength{\evensidemargin}{#2}
                            \setlength{\textwidth}{\paperwidth}
                  \addtolength{\textwidth}{-\oddsidemargin}
                  \addtolength{\textwidth}{-\evensidemargin}
                  \addtolength{\textwidth}{-\marginparsep}
                  \addtolength{\textwidth}{-\marginparwidth}
                            \setlength{\topmargin}{#3}
                            \setlength{\textheight}{\paperheight}
                  \addtolength{\textheight}{-\topmargin}
                  \addtolength{\textheight}{-\headheight}
                  \addtolength{\textheight}{-\headsep}
                  \addtolength{\textheight}{-\footskip}
                  \addtolength{\textheight}{-#4}}
\pageformat{2cm}{3cm}{25mm}{25mm}{1pt}{0pt}

\usepackage{ifthen}
\newboolean{article}
    \setboolean{article}{true}
\newboolean{report}
\newboolean{book}
\newboolean{letter}
\newboolean{german}
\newboolean{italian}
\newboolean{nobaselinestretch}
\newboolean{nosectionappendix}
\newboolean{oldtoc}
\newboolean{nosectionequn}
\newboolean{notheorem}

\ifthenelse{\boolean{german}}{
    \usepackage{german}}{}

\usepackage[latin1]{inputenc}

\usepackage{amsmath}
\usepackage{amssymb}
\usepackage[mathscr]{eucal}

\ifthenelse{\boolean{notheorem}}{}{
    \usepackage{theorem}}

\ifthenelse{\boolean{nobaselinestretch}}{}{
    \renewcommand{\baselinestretch}{1.25}}

\newenvironment{env}[2]{\begin{#1}#2\end{#1}}{}
    \newcommand{\beq}[1]{\begin{env}{equation}{#1}}
    \newcommand{\beqn}[1]{\begin{env}{equation*}{#1}}
    \newcommand{\bal}[1]{\begin{env}{align}{#1}}
    \newcommand{\baln}[1]{\begin{env}{align*}{#1}}
    \newcommand{\bga}[1]{\begin{env}{gather}{#1}}
    \newcommand{\bgan}[1]{\begin{env}{gather*}{#1}}
    \newcommand{\bflal}[1]{\begin{env}{flalign}{#1}}
    \newcommand{\bflaln}[1]{\begin{env}{flalign*}{#1}}
    \newcommand{\bmu}[1]{\begin{env}{multline}{#1}}
    \newcommand{\bmun}[1]{\begin{env}{multline*}{#1}}
    \newcommand{\bsp}[1]{\begin{env}{split}{#1}}

    \newcommand{\eeq}{\end{env}}
    \newcommand{\eeqn}{\end{env}}
    \newcommand{\eal}{\end{env}}
    \newcommand{\ealn}{\end{env}}
    \newcommand{\ega}{\end{env}}
    \newcommand{\egan}{\end{env}}
    \newcommand{\eflal}{\end{env}}
    \newcommand{\eflaln}{\end{env}}
    \newcommand{\emu}{\end{env}}
    \newcommand{\emun}{\end{env}}
    \newcommand{\esp}{\end{env}}

\newcommand{\lf}{\vspace{2ex}}

\renewcommand{\bf}[1]{\textbf{#1}}
\renewcommand{\it}[1]{\textit{#1}}

\renewcommand{\sf}[1]{\textsf{#1}}

\renewcommand{\tt}[1]{\texttt{#1}}
\newcommand{\hl}[1]{\bf{\it{#1}}}

\newcommand{\mbf}[1]{\mathbf{#1}}
\newcommand{\msf}[1]{\text{\small$\sf{#1}$}}
\newcommand{\mtt}[1]{\mathtt{#1}}

\newcommand{\cmc}[1]{\mathcal{#1}}
\newcommand{\eus}[1]{\mathscr{#1}}
\newcommand{\euf}[1]{\mathfrak{#1}}
\newcommand{\bb}[1]{\mathbb{#1}}

\newcommand{\mfootnotesize}[1]{{\setlength{\arraycolsep}{.5ex}\text{\footnotesize$#1$}}}
\newcommand{\mscriptsize}[1]{{\setlength{\arraycolsep}{.3ex}\text{\scriptsize$#1$}}}

\newcommand{\nbd}[1]{$#1$\nobreakdash--}
\newcommand{\ol}[1]{\overline{#1}}

\newcommand{\wh}[1]{\widehat{#1}}

\newcommand{\abs}[1]{\left\lvert#1\right\rvert}

\newcommand{\bfam}[1]{\bigl(#1\bigr)}
\newcommand{\Bfam}[1]{\Bigl(#1\Bigr)}
\newcommand{\AB}[1]{\langle#1\rangle}
\newcommand{\bAB}[1]{\bigl\langle#1\bigr\rangle}

\newcommand{\CB}[1]{\{#1\}}
\newcommand{\bCB}[1]{\bigl\{#1\bigr\}}
\newcommand{\BCB}[1]{\Bigl\{#1\Bigr\}}
\newcommand{\SB}[1]{[#1]}

\newcommand{\Matrix}[1]{\begin{pmatrix}#1\end{pmatrix}}

\newcommand{\fMatrix}[1]{\mfootnotesize{\Matrix{#1}}}
\newcommand{\sMatrix}[1]{\mscriptsize{\Matrix{#1}}}

\newcommand{\set}[2][]{
    \ifthenelse{\equal{#1}{}}{
        \CB{#2}}{
        \CB{#1~|~#2}}}
\newcommand{\bset}[2][]{
    \ifthenelse{\equal{#1}{}}{
        \bCB{#2}}{
        \bCB{#1~|~#2}}}
\newcommand{\Bset}[2][]{
    \ifthenelse{\equal{#1}{}}{
        \BCB{#2}}{
        \BCB{#1~\big|~#2}}}

\DeclareMathOperator{\ls}{\normalfont\msf{span}}

\DeclareMathOperator{\tr}{\normalfont\msf{tr}}
\DeclareMathOperator{\id}{\normalfont\msf{id}}

\renewcommand{\ker}{\operatorname{\msf{ker}}}
\renewcommand{\dim}{\operatorname{\msf{dim}}}
\renewcommand{\Re}{\operatorname{\msf{Re}}}
\renewcommand{\Im}{\operatorname{\msf{Im}}}

\newcommand{\C}{\bb{C}}

\newcommand{\N}{\bb{N}}

\newcommand{\R}{\bb{R}}

\newcommand{\cB}{\cmc{B}}
\newcommand{\cC}{\cmc{C}}
\newcommand{\cD}{\cmc{D}}
\newcommand{\cF}{\cmc{F}}

\newcommand{\sB}{\eus{B}}

\newcommand{\eG}{\euf{G}}
\newcommand{\eH}{\euf{H}}

\newcommand{\U}{\mbf{1}}

\newcommand{\0}{\mbf{0}}

\ifthenelse{\boolean{nosectionequn}}{}{
    \numberwithin{equation}{section}
    }

\ifthenelse{\boolean{article}\or\boolean{letter}\or\boolean{nosectionequn}}{
    \setboolean{nosectionappendix}{true}}{}
\ifthenelse{\boolean{nosectionappendix}}{}{
    \renewcommand{\appendix}{
        \chapter*{\appendixname}
        \addcontentsline{toc}{chapter}{\appendixname}
        \renewcommand{\thesection}{\Alph{section}}
        \setcounter{section}{0}}}
   
\ifthenelse{\boolean{report}\or\boolean{book}}{
    }{}

\ifthenelse{\boolean{notheorem}}{}{
        \newcommand{\mnname}{Mathematical note.}
        \newcommand{\enname}{End of the note.}
        \newcommand{\definame}{Definition.}
        \newcommand{\propname}{Proposition.}
        \newcommand{\lemname}{Lemma.}
        \newcommand{\exname}{Example.}
        \newcommand{\exername}{Exercise.}
        \newcommand{\remname}{Remark.}
        \newcommand{\obname}{Observation.}
        \newcommand{\thmname}{Theorem.}
        \newcommand{\corname}{Corollary.}
        \newcommand{\proofname}{Proof.}
    \ifthenelse{\boolean{german}}{
        \renewcommand{\mnname}{Mathematische Notiz.}
        \renewcommand{\enname}{Ende der Notiz.}
        \renewcommand{\exname}{Beispiel.}
        \renewcommand{\exername}{Übung.}
        \renewcommand{\remname}{Bemerkung.}
        \renewcommand{\obname}{Beobachtung.}
        \renewcommand{\thmname}{Satz.}
        \renewcommand{\corname}{Korollar.}
        \renewcommand{\proofname}{Beweis.}}{}
    \ifthenelse{\boolean{italian}}{
        \renewcommand{\mnname}{Nota matematica.}
        \renewcommand{\enname}{Fina della nota.}
        \renewcommand{\definame}{Definizione.}
        \renewcommand{\propname}{Proposizione.}
        \renewcommand{\exname}{Esempio.}
        \renewcommand{\exername}{Esercizio.}
        \renewcommand{\remname}{Nota.}
        \renewcommand{\obname}{Osservazione.}
        \renewcommand{\thmname}{Teorema.}
        \renewcommand{\corname}{Corollario.}
        \renewcommand{\proofname}{Dimostrazione.}

       \renewcommand{\appendixname}{Appendice}

       }{}
    \theoremheaderfont{\normalfont\bfseries}
    \theoremstyle{change}
        \theorembodyfont{\rmfamily}
            \newtheorem{emp}{}[section]
                \newcommand{\bemp}[1][]{
                    \begin{emp}\hskip-\labelsep\bf{#1}\hskip\labelsep}
                \newcommand{\eemp}{\end{emp}}
\newtheorem{itemp}[emp]{}
                \newcommand{\bitemp}[1][]{
                    \begin{itemp}\hskip-\labelsep\bf{#1}\hskip\labelsep\normalfont\itshape}
                \newcommand{\eitemp}{\end{itemp}}
            \newtheorem{mn}[emp]{\mnname}
                \newcommand{\bnm}{\begin{mn}~\begin{quotation}\renewcommand{\baselinestretch}{1}\small\noindent\ignorespaces}
                \newcommand{\enm}{\end{quotation}\hfill\bf{\enname}\end{mn}}
            \newtheorem{ex}[emp]{\exname}
                \newcommand{\bex}{\begin{ex}}
                \newcommand{\eex}{\end{ex}}
            \newtheorem{exer}[emp]{\exername}
                \newcommand{\bexer}{\begin{exer}}
                \newcommand{\eexer}{\end{exer}}
            \newtheorem{defi}[emp]{\definame}
                \newcommand{\bdefi}{\begin{defi}}
                \newcommand{\edefi}{\end{defi}}
            \newtheorem{rem}[emp]{\remname}
                \newcommand{\brem}{\begin{rem}}
                \newcommand{\erem}{\end{rem}}
            \newtheorem{ob}[emp]{\obname}
                \newcommand{\bob}{\begin{ob}}
                \newcommand{\eob}{\end{ob}}
        \theorembodyfont{\normalfont\itshape}
            \newtheorem{thm}[emp]{\thmname}
                \newcommand{\bthm}{\begin{thm}}
                \newcommand{\ethm}{\end{thm}}
            \newtheorem{prop}[emp]{\propname}
                \newcommand{\bprop}{\begin{prop}}
                \newcommand{\eprop}{\end{prop}}
            \newtheorem{cor}[emp]{\corname}
                \newcommand{\bcor}{\begin{cor}}
                \newcommand{\ecor}{\end{cor}}
            \newtheorem{lem}[emp]{\lemname}
                \newcommand{\blem}{\begin{lem}}
                \newcommand{\elem}{\end{lem}}
\newenvironment{empn}[1]{\lf\noindent\bf{#1}\ignorespaces\hskip\labelsep}{\lf}
        \newcommand{\bempn}[1]{\begin{empn}{#1}}
        \newcommand{\eempn}{\end{empn}}
        \newcommand{\bitempn}[1]{\begin{empn}{#1}\normalfont\itshape}
        \newcommand{\eitempn}{\end{empn}}
                \newcommand{\bnmn}{\begin{empn}{\mnname}~\begin{quotation}\renewcommand{\baselinestretch}{1}\small\noindent\ignorespaces}
                \newcommand{\enmn}{\end{quotation}\hfill\bf{\enname}\end{empn}}
        \newcommand{\bexn}{\begin{empn}{\exname}}
        \newcommand{\eexn}{\end{empn}}
        \newcommand{\bexern}{\begin{empn}{\exername}}
        \newcommand{\eexern}{\end{empn}}
        \newcommand{\bdefin}{\begin{empn}{\definame}}
        \newcommand{\edefin}{\end{empn}}
        \newcommand{\bremn}{\begin{empn}{\remname}}
        \newcommand{\eremn}{\end{empn}}
        \newcommand{\bobn}{\begin{empn}{\obname}}
        \newcommand{\eobn}{\end{empn}}

\newcommand{\qedsymbol}{~\rule[-0.35mm]{2mm}{2mm}}
    \newcounter{proof}[emp]
    \newenvironment{Proof}[1]{
        \vspace{1ex}
        \renewcommand{\item}[1][\stepcounter{proof}(\roman{proof})]%
            {##1\hskip\labelsep}
        \noindent\textsc{#1\hskip\labelsep}}{
        \nolinebreak\qedsymbol}
    \newcommand{\proof}[1][\proofname]{
        \begin{Proof}{#1}\ignorespaces}
    \newcommand{\qed}{\end{Proof}}
    \newcommand{\noqed}{
        \renewcommand{\qedsymbol}{}
        \end{Proof}}}
    \ifthenelse{\boolean{italian}}{
        \renewcommand{\proofname}{Dimostrazione.}}{}

\usepackage[varg]{txfonts}

\usepackage[hypertex]{hyperref}

\begin{document}

\title{Maximal Commutative Subalgebras Invariant for CP-Maps: (Counter-)Examples\thanks{BVRB is supported by INDAM (Italy) and UKIERI (UK). FF and MS are supported by the Italian MUR (PRIN 2005) and by GNAMPA (``Semigruppi Markoviani Quantistici'' 2008). MS is supported by research funds of the Dipartimento S.E.G.e S. of University of Molise.}}

\author{B.V.\ Rajarama Bhat, Franco Fagnola, and Michael Skeide}

\date{April 2008}

\maketitle

\vspace{-6ex}
\begin{abstract}
\noindent
We solve, mainly by counterexamples, many natural questions regarding maximal commutative subalgebras invariant under CP-maps or semigroups of CP-maps on a von Neumann algebra. In particular, we discuss the structure of the generators of norm continuous semigroups on $\sB(G)$ leaving a maximal commutative subalgebra invariant and show that there exists Markov CP-semigroups on $M_d$ without invariant maximal commutative subalgebras for any $d>2$.\vspace{-3ex}
\end{abstract}

\section{Introduction}

\hl{Markov semigroups}, that is, semigroups of normal unital completely positive (CP-)maps on a von Neumann algebra $\cB\subset\sB(G)$ ($G$ a  Hilbert space) are models for irreversible evolutions both of classical and of quantum systems. 
Indeed, if $G$ is separable, then a commutative von Neumann algebra $\cC\subset\sB(G)$ is isomorphic to $L^\infty(\Omega,\cF,P)$ for some probability space, and a Markov semigroup on $\cC$ is the semigroup induced on $L^\infty(\Omega,\cF,P)$ by a classical Markov semigroup of transition probabilities. More generally, if a Markov semigroup $T=\bfam{T_t}_{t\in\R_+}$ on a not necessarily commutative von Neumann algebra $\cB$ leaves a commutative subalgebra $\cC$ \hl{invariant} (that is, $T_t(\cC)\subset\cC$ for all $t\in\R_+$), then the restriction to $\cC$ gives rise to a classical Markov semigroup. Finding invariant commutative subalgebras means, thus, recognizing classical subsystems as embedded into a quantum one.

The study of invariant commutative subalgebras initiated in 1989 in the framework of quantum flows when P.-A. Meyer wrote the short note \cite{Mey89p} showing how certain finite Markov chains in continuous time can be expressed as quantum flow in Fock space. Meyer's construction was extended by Parthasarathy and Sinha in \cite{PaSi90} by constructing the structure maps of the flow through certain group actions. Later Fagnola showed (see, e.g., \cite{Fag99}) that also classical diffusion processes can be viewed as restrictions to a commutative subalgebras of a quantum flow. Quantum Markov flows and semigroups with an invariant commutative subalgebra (the algebra generated by the system Hamiltonian) arise in a natural way in the stochastic limit; many examples can be found in the book \cite{ALV02} by Accardi, Lu and Volovich.

The above investigations, either by construction or as a result of a scaling limit of a Hamiltonian evolution, lead to a quantum Markov flow (respectively semigroup) on a $\sB(G)$ with a restriction to a commutative subalgebra $\cC$ coinciding with the flow (respectively semigroup) of a prescribed classical Markov process. The more difficult problem of characterizing \emph{all} the invariant commutative subalgebras of a given quantum flow (respectively semigroup), however, was not attacked.

Recently, Rebolledo \cite{Reb05}, motivated by the interpretation of decoherence as the appearance of classical features in quantum evolutions, found a simple sufficient algebraic condition for finding a maximal abelian subalgebra invariant under the action of a quantum Markov semigroup.

This paper is concerned with the problem of finding all invariant maximal commutative subalgebras $\cC$ of a CP-semigroup on $\cB\subset\sB(G)$ and of its generator. 

A commutative subalgebra $\cC$ with $\U=\U_\cB\in\cC\subset\cB\subset\sB(G)$ is a \hl{maximal commutative subalgebra} of $\cB$, if $\cC\subset\cD\subset\cB$ for a commutative subalgebra $\cD$ implies $\cD=\cC$. A maximal commutative subalgebra of $\cB=\sB(G)$ is a called a \hl{maximal abelian subalgebra} or a \hl{masa}. If $G$ is separable and if $\cC\subset\sB(G)$ is a masa isomorphic to $L^\infty(\Omega,\cF,P)$, then $G\cong L^2(\Omega,\cF,P)$. If $\cC$ is a maximal comutative subalgebra of $\cB\subset\sB(G)$, then we obtain a description of the system by \it{classical} (or \it{macroscopic}) parameters that is not improvable by measuring a set of classical observables. If $\cC$ is a masa, then this description is complete.

Rebolledo \cite{Reb05} (see also \cite{Reb05a}) proved the following sufficient criterion in the case $\cB=\sB(G)$: Let $T$ be a normal CP-map on $\sB(G)$ given by some \hl{Kraus decomposition} $T(b)=\sum_iL_i^*bL_i$ $(L_i\in\sB(G))$. Suppose that $\cC\subset\sB(G)$ is a masa generated by a single self-adjoint element $c\in\cC$, and suppose that there are self-adjoint elements $c_i\in\cC$ such that
\beqn{
cL_i-L_ic
~=~
c_iL_i.
}\eeqn
Then $T(\cC)\subset\cC$. If $T$ is the \hl{CP-part} of the \hl{Gorini-Kossakowski-Sudarshan-Lindblad generator} \cite{GKS76,Lin76}
\beqn{
L(b)
~=~
\sum_iL_i^*bL_i+b\beta+\beta^*b
}\eeqn
$(\beta\in\sB(G))$ of a uniformly continuous CP-semigroup $T_t=e^{tL}$ on $\sB(G)$, then invariance of the CP-part plus invariance of the \hl{effective Hamitonian} $b\mapsto b\beta+\beta^*b$ implies that the whole CP-semigroup leaves $\cC$ invariant. In the case of a Markov semigroup (where $L$ has to be normalized to $L(\U)=0$) we get
\beqn{
L(b)
~=~
\sum_iL_i^*bL_i-{\textstyle\frac{b\bfam{\sum_iL_i^*L_i}+\bfam{\sum_iL_i^*L_i}b}{2}}+i\SB{b,h},
}\eeqn
for the self-adjoint $h=\Im\beta\in\sB(G)$. As the CP-part $T$ alone, by Rebolledo's criterion, leaves $\cC$ invariant, we have, in particular, that $\sum_iL_i^*L_i=T(\U)\in\cC$. So, if (and only if; see \cite[Lemma 4.4]{FaSk07}) also $h\in\cC$ so that the \hl{Hamiltonian} $b\mapsto i\SB{b,h}$ leaves $\cC$ invariant, then all $T_t$ leave $\cC$ invariant.

Fagnola and Skeide \cite{FaSk07} proved the following generalization of Rebolledo, which now provides a sufficient and necessary criterion.

\bitemp[Theorem \cite{FaSk07}.~]\label{Rebcor}
Let $T$ be a normal CP-map on $\sB(G)$ with Kraus decomposition $T(b)=\sum_{i\in I}L_i^*bL_i$. Then $T$ leaves a maximal abelian von Neumann algebra $\cC\subset\sB(G)$ invariant, if and only if for every $c\in\cC$ there exist coefficients $c_{ij}(c)\in\cC$ $(i,j\in I)$ such that
\baln{
1.)~~~
&
c_{ij}(c^*)
~=~
c_{ji}(c)^*,
&
2.)~~~
&
cL_i-L_ic
~=~
\sum_{j\in I}c_{ij}(c)L_j,
}\ealn
for all $c\in\cC$.
\eitemp

Theorem \ref{Rebcor} is a special case of \cite[Theorem 3.1]{FaSk07} for general von Neumann algebras. Fagnola and Skeide also provide the sufficient and necessary criterion \cite[Theorem 4.2]{FaSk07} for the generator of a uniformly continuous CP-semigroup on a general von Neumann algebra. We state here the result of the specialization to $\sB(G)$. A proof is delegated to the appendix.

\bthm\label{Lincor}
Let $L$ be the generator of a uniformly continuous normal CP-semigroup on $\sB(G)$ with Gorini-Kossakowski-Sudarshan-Lindblad form $L(b)=\sum_{i\in I}L_i^*bL_i+b\beta+\beta^* b$. Then $L$, or equivalently, all $T_t=e^{tL}$, leave a maximal abelian von Neumann algebra $\cC\subset\sB(G)$ invariant, if and only if there exist coefficients $\gamma=\gamma^*,c_i\in\cC$, and for every $c\in\cC$ there exist coefficients $c_{ij}(c)\in\cC$ $(i,j\in I)$ such that
\baln{
1.)~~~
&
c_{ij}(c^*)
~=~
c_{ji}(c)^*,
&
2.)~~~
&
cL_i-L_ic
~=~
\sum_{j\in I}c_{ij}(c)(L_j-c_j),
}\ealn
\vspace{-4ex}
\baln{
3.)~~~
&
L(c)
~=~
\sum_{i\in I}(L_i-c_i)^*c(L_i-c_i)+\gamma c
}\ealn
for all $c\in\cC$.
\ethm

\brem
We would like to mention that in both theorems (like in  Theorems \ref{TcharCCPthm} and \ref{LcharCCPthm}, from which the former are derived) maximal commutativity of $\cC$ easily guarantees sufficiency. The stated conditions are necessary (in all four theorems) for invariance of the unital commutative subalgebra $\cC$, even if $\cC$ is not maximal commutative.
\erem

Like \cite[Theorem 30.16]{Par92}, the following theorem characterizes the possibilities to transform a generator in minimal Gorini-Kossakowski-Sudarshan-Lindblad form into another. The proof illustrates the power of techniques from product systems of Hilbert modules. But as we do not need these techniques in the rest of these notes, we postpone also this proof to the appendix.

\bthm\label{Linchthm}
Let $L$ be the generator of a uniformly continuous normal CP-semigroup on $\sB(G)$ in minimal Gorini-Kossakowski-Sudarshan-Lindblad form $L(b)=\sum_{i\in I}L_i^*bL_i+b\beta+\beta^* b$, and let $K(b)=\sum_{j\in J}K_j^*bK_j+b\alpha+\alpha^* b$ be another generator.

Then $K=L$, if and only if there exists a matrix $\sMatrix{\gamma&\eta^*\\\eta'&M}\in M_{(1+\#J)\times(1+\#I)}$, with $\eta'\in\C^{\#J}$ arbitrary, $M=\bfam{a_{ji}}_{ji}\in M_{\#J\times\#I}$ an isometry, $\eta=-M^*\eta'\in\C^{\#I}$, and $\gamma=ih-\frac{\AB{\eta',\eta'}}{2}\in\C$ ($h\in\R$ arbitrary), such that
\baln{
\alpha
&
~=~
\beta+\gamma\U+\sum_{i\in I}\ol{\eta}_iL_i,
&
K_j
&
~=~
\eta'_j\U+\sum_{i\in I}a_{ji}L_i.
}\ealn
This holds for arbitrary cardinalities $\#I$ and $\#J$, if infinite sums are understood as strongly convergent.
\ethm

\bcor
A similar result holds if the Gorini-Kossakowski-Sudarshan-Lindblad form of $L$ is not necessarily minimal. In that case $M$ may be just a partial isometry and $\eta'$ must be such that $MM^*\eta'=\eta'$.
\ecor

\proof
Observe that the minimal $L_i$ in the theorem may be recoverd as $L_i=\eta_i\U+\sum_{j\in J}\ol{a}_{ji}K_j$. So, in order to compare two not necessarily minimal Gorini-Kossakowski-Sudarshan-Lindblad forms we may simply ``factor'' through a minimal one.\qed

\lf
There are several natural questions around about Theorems \ref{Rebcor} and \ref{Lincor} and how they are related with Rebolledo's original criterion. Most of them are motivated by the examples with \nbd{2\times 2}matrices that have been studied in \cite{FaSk07}. The goal of these notes is to give answers to these questions, and Theorem \ref{Linchthm} will play a crucial role. As our results here show, the answers sometimes are typical only for $M_2$ and look different already for $M_3$. Therefore, in the following list of questions and throughout the answers later on in these notes we will have to distinguish between $M_2$ and higher dimensional settings.

We explain briefly why counterexamples for a single map furnish also counterexamples for the semigroup case.

\bob\label{mapSGob}
The CP-semigroup $T_t=e^{tL}$ leaves a subalgebra invariant, if and only if its generator $L$ leaves that subalgebra invariant. So, for all questions about invariance for CP-semigroups we are done if we answer the single mapping case. (If $T$ is CP-map with a certain invariance property, then $e^{tT}$ shares that property.) Similarly, if $T$ is a \bf{unital} CP-map leaving a certain subalgebra invariant or not, then $L:=T-\id$ is the generator of a Markov semigroup sharing this property. (This is so, simply because $\id$ leaves every subalgebra invariant so that $L$ and, therefore, $T_t$ share the invariance properties of $T$.)
\eob

We do not know, whether the converse statement of Observation \ref{mapSGob} is also true. (If $L$ leaves no masa invariant invariant, then Observation \ref{mapSGob} tells us that no masa is invariant for \bf{all} $T_t$. But \it{a priori} it might be possible that $T$ has ``wandering'' invariant masas.)

\lf
We now list our questions and the answers we obtain later on in the remainder of these notes.

\begin{enumerate}
\item\label{Q1}
Does every CP-semigroup on $\sB(G)$ leave some masa invariant?

Answer: No, by Example \ref{M2noninvex} already for a single CP-map on $M_2$ and, therefore, also for a CP-semigroup on $M_2$ (and, therefore, for all $\sB(G)$).

\item\label{Q2}
Does every Markov semigroup on $\sB(G)$ leave some masa invariant?

Answer: Yes, for $M_2$ by Theorem \ref{M2invthm} both for Markov semigroups and for unital CP-maps.

Answer: No, for $M_d$ $(d\ge3)$ and for $\sB(G)$ by Example \ref{M3noninvex} for Markov semigroups, even in countably infinite dimension, and no for unital CP-maps by Example \ref{M3noninvex'} in finite dimension.

\item\label{Q3}
Is Rebolledo's criterion equivalent to the one in Theorem \ref{Rebcor}? More precisely, does every normal CP-map on $\sB(G)$ that leaves a masa invariant, admit a Kraus decomposition fulfilling Rebolledo's criterion?

Answer: No, already for $M_2$, by Example \ref{nonRebex} for Markov semigroups and for unital CP-maps.

\item\label{Q4}
Suppose we have a generator leaving a masa invariant. Does every such generator decompose, like in Rebolledo's criterion, into a CP-part and a Hamiltonian part that leave the masa invariant, separately?

Answer: No for the CP-part, already in the case of a Markov semigroup on $M_2$ by Example \ref{CP2nonex}. This answer extends to all $\sB(G)$.

Answer: Yes for the Hamiltonian part, in the case of CP-semigroups on $M_2$ by Corollary \ref{H2cor}. No, in the case of CP-semigroups on $M_3$ and higher dimension, by Example \ref{H3nonex}.
\end{enumerate}

\noindent
In Section \ref{M2SEC} we study everything related to $\cB=M_2$, while Section \ref{M3SEC} is dedicated to $\cB=M_d$ $(d\ge3)$ and the infinite-dimensional case.

We would like to mention that a further natural question asked in \cite{FaSk07}, namely, whether the necessary and sufficient criterion in \cite{FaSk07} remains valid for unbounded generators, has a negative answer, too. There exist generators in terms of double commutators and the CCR that leave invariant a masa but that do not fufill the (unbounded analogue of the) criterion in \cite{FaSk07}. We will study these generators  elsewhere systematically. Here we restrict ourselves to the bounded case.

We also mention also that the relationships we find in Theorem and \ref{Lincor}, Theorem \ref{Linchthm} and its corollary among the operators appearing in the Gorini-Kossakowski-Sudarshan-Lindblad representation of a generator are new. These, together with those satisfied by generators of other special classes of quantum Markov semigroups (see, e.g., \cite{Dav79,Hol96,BaPa96,AlGo02,AFH06,FaUm07}), reveal the rich algebraic structure of generators of CP-semigroups.  

\lf\noindent
\bf{Conventions.~} For every $n\in\N$ we denote by $M_n=M_n(\C)=\sB(\C^n)$ the von Neumann algebra of \nbd{n\times n}matrices with complex entries. By $M_\infty$ we mean the von Neumann algebra $\sB(G)$ for a separable infinite-dimensional Hilbert space $G$. The elements of $\sB(G)$ are considered as matrices with respect to a fixed orthonormal basis $\bfam{e_n}_{n\in\N_0}$ of $G$. By $\cD_n$ $(n\in\N\cup\CB{\infty})$ we denote the respective subalgebras of diagonal matrices.

\lf\noindent
\bf{Acknowledgements.~}
BVRB is grateful to  FF, MS and L.\ Accardi for their generous hospitality during his visit to Italy in September-October 2008. Most of Section \ref{M3SEC} was done during that visit. MS wishes to thank BVRB for a nice stay at ISI Bangalore during December 2007 to February 2008, where the first part of these notes has been written.

\section{Examples and results for $M_2$}\label{M2SEC}

We start  with some counterexamples for things that do not even work for $M_2$.

\bex\label{M2noninvex}
Consider the CP map $T\colon M_2\rightarrow M_2$ defined by
\beqn{
T\fMatrix{
a & b
\\
c & d
}
~=~
\fMatrix{
1 & 1 
\\
0 & 1
}
\fMatrix{
a & b 
\\
c & d
}
\fMatrix{
1 & 0 
\\
1 & 1
}
~=~
\fMatrix{
a+b+c+d & b+d 
\\
c+d & d
}
.
}\eeqn
If $T$ leaves a masa $\cC\subset M_2$ invariant, then $\CB{\U,T(\U),T^2(\U),\ldots}\subset\cC$ should all commute. But clearly this is not the case as $T(\U)=\sMatrix{2 & 1 \\1 & 1}$ and $T^2(\U)=\sMatrix{5 & 2 \\2 & 1}$ do not commute. So neither the CP-map $T$ nor the CP-semigroup $e^{tT}$ leave a masa of $M_2$ invariant.

Any CP-map $T$ may be extended to a CP-map $\wh{T}(X)=T(\U_2 X\U_2)$ on $M_d$ for any $d\ge3$ including $\infty$. Again $\wh{T}(\U)$ and $\wh{T}^2(\U)$ do not commute. So, $\wh{T}$ has no invariant masa and the CP-semigroup $e^{t\wh{T}}$ $(=\wh{e^{tT}})$ shares this property.
\eex

\bex\label{nonRebex}
Define the CP-map $T\colon M_2\rightarrow M_2$ by
\beqn{
T(X)
~=~
L_1^* X L_1 + L_2^\ast X L_2
}\eeqn
where
\baln{
L_1
&
~=~
\Matrix{
\frac{1}{\sqrt{2}} & 0 
\\
\frac{1}{2} & \frac{1}{2}
}
,
&
L_2
&
~=~
\Matrix{
0 & \frac{1}{\sqrt{2}} 
\\
-\frac{1}{2} & \frac{1}{2}
}
.
}\ealn
Then
\beqn{
T
\Matrix{
a & b 
\\
c & d 
}
~=~
\Matrix{
\frac{1}{2}(a+\frac{b}{\sqrt{2}}+\frac{c}{\sqrt {2}}+d) & \frac{b-c}{2 \sqrt{2}} 
\\
\frac{c-b}{2 \sqrt{2}} & \frac{1}{2}(a+\frac{b}{\sqrt{2}}+\frac{c}{\sqrt {2}}+d)
}
.
}\eeqn
We see that $T$ is unital and that it leaves the diagonal subalgebra $\cD_2$ of $M_2$ invariant.

Now suppose $T(X) = \sum  \limits_j
K^*_j X K_j$ is another Kraus decomposition of $T$. Then
each $K_j,$ is a linear combination of $L_1, L_2$; see Observation \ref{lincob}. Say $K_1 = aL_1
+bL_2,~~a,b \in \C.$  Now suppose this decomposition
satisfies Rebolledo's condition.  Then for every diagonal matrix
$D=\sMatrix{d_1 & 0 \\ 0 & d_2}\in\cD_2$ there exists $D' = \sMatrix{d'_1 & 0 \\ 0 & d'_2}\in\cD_2$ (depending upon $D$) such that
\beqn{
D'K_1
~=~
K_1D.
}\eeqn
So
\beqn{
\Matrix{
d'_1 & 0 
\\
0 & d'_2
}
\Matrix{
\frac{a}{\sqrt{2}} & \frac{b}{\sqrt{2}} 
\\
\frac{a-b}{2} & \frac{a+b}{2} }
~=~
\Matrix{
\frac{a}{\sqrt{2}} & \frac{b}{\sqrt{2}} 
\\
\frac{a-b}{2} & \frac{a+b}{2} 
}
\Matrix{
d_1 & 0 
\\
0 & d_2
}
}\eeqn
or
\beqn{
\Matrix{
d'_1 \frac{a}{\sqrt{2}} & d'_1 \frac{b}{\sqrt{2}} 
\\
d'_2\left(\frac{a-b}{2}\right) & d'_2\left(\frac{a+b}{2}\right)
}
~=~
\Matrix{
d_1 ~\frac{a}{\sqrt{2}} & d_2~ \frac{b}{\sqrt{2}} 
\\
d_1 \left(\frac{a-b}{2}\right) & d_2 \left(\frac{a+b}{2}\right)
}
}\eeqn
It is easily seen that no non-zero $K_1$ satisfies this condition. We conclude that Rebolledo's condition is not a necessary condition.
\eex

We now discuss several things that work only for $M_2$. The counterexamples in the general case for the statements we prove here for $M_2$, must wait until Section \ref{M3SEC} on $M_3$.

\blem
Let $\alpha$ be a linear \nbd{*}map on $M_2$ such that $\alpha(\U)\in\C\U$. Then $\alpha$ leaves a masa of $M_2$ invariant.
\elem

\proof
The \it{Cayley-Hamilton theorem} asserts that for every matrix $Y\in M_n$ the characteristic polynomial $P$ of $Y$ gives $P(Y)=0$. It follows that for every $Y\in M_2$ the subalgebra of $M_2$ generated by $Y$ has the form $\cC_Y:=\C\U_2+\C Y$. Therefore, if we find a self-adjoint $Y=Y^*\notin\C\U_2$ such that $\alpha(Y)\in\cC_Y$, then $\cC_Y$ is a masa of $M_2$ invariant for $\alpha$.

Define the \nbd{4}dimensional real subspace $S=\CB{X\in M_2\colon X=X^*}$ of self-adjoint elements of $M_2$. By $\tr$ we denote the normalized trace on $M_2$. Then $\id_S-\tr\U\colon X\mapsto X-\tr(X)\U$ defines a projection onto the subspace $S_0:=S\cap\ker\,\tr$ of self-adjoint zero-trace operators. The linear map $(\id_S-\tr\U)\circ\alpha$ leaves the \nbd{3}dimensional real vector space $S_0$ invariant. Therefore, $\beta:=(\id_S-\tr\U)\circ\alpha\upharpoonright S_0$ has an eigenvector $Y$ to some real eigenvalue. Clearly, $\alpha(Y)\in\cC_Y$ and $Y\notin\C\U_2$, so that $\cC_Y$ is a masa invariant for $\alpha$. (Of course, it is an easy exercise to check directly that $Y^2\in\C\U$ for every self-adjoint zero-trace operator $Y\in M_2$, showing that $\cC_Y$ is an algebra without reference to the Cayley-Hamilton theorem.)\qed

\lf
The following theorem is a simple corollary of the lemma.

\bthm\label{M2invthm}
Every unital CP-map $T$ on $M_2$ has an invariant masa. Every generator $L$ of a Markov semigroup on $M_2$ has an invariant masa. 
\ethm

\proof
$T$ is a linear \nbd{*}map that maps $\U$ to $1\cdot\U$ and $L$ is a linear \nbd{*}map that maps $\U$ to $0\cdot\U$.\qed

\lf
Once assured existence of an invariant masa of $M_2$, by a basis transformation we may always assume that this invariant subalgebra is $\cD_2$. We now investigate when a generator leaving $\cD_2$ invariant can be split such that also its CP-part or at least its Hamiltonian part leaves $\cD_2$ invariant. Note that by Corollary \ref{H2cor} and Example \ref{CP2nonex} these two properties need not coincide.

\bthm\label{D2invthm}
Suppose the minimal Gorini-Kossakowski-Sudarshan-Lindblad generator $L(X)$ $=\sum_{i=1}^dL_i^*XL_i+XB+B^*X$ of a CP-semigroup on $M_n$ leaves $\cD_n$ invariant. Then $L$ admits a (minimal) Gorini-Kossakowski-Sudarshan-Lindblad form whose CP-part leaves $\cD_n$ invariant separately, if and only if there is a linear combination $K:=\sum_{i=1}^d\ol{\eta}_iL_i$ such that $B+K\in\cD_n$.
\ethm

\proof
Note that $T(X)=\sum_{i=1}^dL_i^*XL_i$ leaves $\cD_n$ invariant, if and  only if $\Delta\colon X\mapsto XB+B^*X=\CB{X,\Re B}+i\SB{X,\Im B}$ does. We show that this happens, if and only if $B\in\cD_n$. ``If'' being clear, for ``only if'' suppose that $\Delta$ leaves $\cD_n$ invariant. Then $\Delta(\U)=2\Re B\in\cD_n$ and, therefore $\CB{X,\Re B}\in\cD_n$ for all $X\in\cD_n$. Clearly, $X\mapsto\SB{X,\Im B}$ leaves $\cD_n$ invariant, if and only if $\Im B\in\cD_n$.

Suppose $A,K_j$ are the coefficients of another Gorini-Kossakowski-Sudarshan-Lindblad form of $L$. So, in order that the CP-part $\sum_{j=1}^{d'}K_j^*XK_j$ leaves $\cD_n$ invariant, it is necessary and sufficient that $A\in\cD_n$. By Theorem \ref{Linchthm} the only possibility to achieve this, is adding linear combinations of the $L_i$ (and $\U$) to $B$. So, the condition $\exists K=\sum_{i=1}^d\ol{\eta}_iL_i\colon B+K\in\cD_n$ is necessary. On the other hand, suppose that $K$ exists. In view of Theorem \ref{Linchthm} put $M=\U_d$, $\eta'=-\eta$, and $\gamma=-\frac{\AB{\eta,\eta}}{2}$. Then the Gorini-Kossakowski-Sudarshan-Lindblad generator with coefficients $A=B+\gamma\U+K$ and $K_i=L_i-\eta_i\U$ coincides with $L$ and $X\mapsto XA+A^*X$ leaves $\cD_n$ invariant.\qed

\bcor\label{H2cor}
Every generator $L$ of a CP-semigroup on $M_2$ leaving $\cD_2$ invariant can be written in a Gorini-Kossakowski-Sudarshan-Lindblad form where also the Hamiltonian part leaves $\cD_2$ invariant.
\ecor

\proof
Either all $L_i$ are in $\cD_2$ so that also $B\in\cD_2$, or there exists at least one $L_k$ that is not diagonal. $H=\Im B$ is self-adjoint and $\Im(\ol{\eta}_kL_k)$ will eliminate the off-diagonal from $H$ for suitable $\eta_k$.\qed

\brem
Note that this is true for arbitrary generators (not necessarily leaving $\cD_2$ invariant) as soon as the CP-part does not leave $\cD_2$ invariant (assuring existence of a nondiagonal $L_k$).
\erem

\bex\label{CP2nonex}
Let $L(X)=L_1^*XL_1+L_2^*XL_2+XB+B^*X$ with
\baln{
B
&
~:=~
-\textstyle\frac{1}{2}\fMatrix{7&6\\10&8},
&
L_1
&
~:=~
\fMatrix{1&1\\1&1},
&
L_2
&
~:=~
\fMatrix{1&2\\2&2}.
}\ealn
One easily verifies that $L$ leaves $\cD_2$ invariant and that $L(\U)=0$. However, all linear combinations of $L_1$ and $L_2$ have equal off-diagonal elements, and $B$ has not. Therefore, none of the linear combinations $B+\gamma\U+\ol{\eta}_1L_1+ \ol{\eta}_2L_2$ will be diagonal. In conclusion, it is not possible to find a Gorini-Kossakowski-Sudarshan-Lindblad form with effective Hamiltonian and CP-part that leave $\cD_2$ invariant separately.

This example extends easily to arbitrary higher dimension $\sB(G)$, if we embed all coefficients it into the \nbd{M_2}corner of $\sB(G)$.
\eex

\section{Examples for $d\ge3$}\label{M3SEC}

Apart from the counterexamples, the preceding section contained also some positive results which were, however, specific for $M_2$. In the present section we give counterexamples to the analogue statements in $M_3$.

\brem
This behaviour, a qualitative jump for what is possible when passing from dimension $2$ to dimension $3$, though not untypical, provided us with some surprises. (In \cite{FaSk07} only two-dimensional examples were studied.) As the comparably ``large'' numbers in Example \ref{H3nonex} and and the more sophisticated construction of Examples \ref{M3noninvex} and \ref{M3noninvex'} show, that these examples were not exactly obvious.
\erem

We start with an example in $M_3$ that contradicts the statement of Corollary \ref{H2cor} for $M_2$.

\bex\label{H3nonex}
Let $L(X)=L_1^*XL_1+L_2^*XL_2+XB+B^*X$ with
\baln{
B
&
~:=~
-\frac{1}{2}\fMatrix{7&6&0\\2&11&0\\4&10&26},
&
L_1
&
~:=~
\fMatrix{1&3&0\\1&0&0\\0&1&5},
&
L_2
&
~:=~
\fMatrix{0&0&0\\1&1&0\\2&0&1}.
}\ealn
One calculates
\beqn{
L\fMatrix{d_1&0&0\\0&d_2&0\\0&0&d_3}
~=~
\fMatrix{-6d_1+2d_2+4d_3&0&0\\0&~9d_1-10d_2+d_3~&0\\0&0&0},
}\eeqn
so that $L$ leaves $\cD_3$ invariant and $L(\U)=0$. One easily computes
\baln{
2B'
~:=~
c_1L_1+c_2L_2+2B
~=~
\fMatrix{c_1-7&3c_1-6&0\\c_1+c_2-2&c_2-11&0\\2c_2-4&~c_1-10~&~5c_1+c_2-26}.
}\ealn
For that $B'-B'^*$ is diagonal we obtain the three equations $c_1-10=0$, $2c_2-4=0$, and $3c_1-6-\ol{c}_1-\ol{c}_2+2=0$. Inserting $c_1=10$ and $c_2=2$ into the third equation gives $30-6-10-2+2=14\ne0$. We conclude that no other Gorini-Kossakowski-Sudarshan-Lindblad form of $L$ has a Hamiltonian part leaving $\cD_3$ invariant.

Also this example extends easily to arbitrary higher dimensional $\sB(G)$, if we embed all coefficients into the \nbd{M_3}corner of $\sB(G)$.
\eex

We now construct examples in $M_d$ $(d\ge3)$ that contradict the statements of Theorem \ref{M2invthm} for $M_2$. The first example is for Markov semigroups and works also for $d=\infty$ (separable!). The second example works for unital CP-maps and, therefore, also for Markov semigroups, but, so far, only for finite $d\ge3$. The idea, common to both examples, is to start with a CP-part that has the simplest possible structure: Compression with a rank-one projection. Then, perturb it with a \it{Hamiltonian perturbation}. If the Hamiltonian has the worst commuting behaviour possible with the rank-one projection, then we obtain a counter example.

\bex\label{M3noninvex}
Let $G$ be a Hilbert space of dimension $d\in\N\cup\CB{\infty}$. Choose a unit vector $e$ and a self-adjoint element $H\in\sB(G)$ such that $\CB{H}'\cap\CB{ee^*}'=\C\U$. (Note that this means that $G$ is separable. Indeed, $G_{H,e}^\perp:=\CB{H^ne\colon n\in\N_0}^\perp$ is a subspace invariant for $H$, and nonseparable if $G$ is. If $H$ generates a masa of $\sB(G)$, then the restriction of $H$ to $G_{H,e}^\perp$ generates a masa of $\sB(G_{H,e}^\perp)$. So, $G_{H,e}^\perp$ and, therefore, also $G$ cannot be nonseparable.) It is easy to write down concrete $G,e,H$ fulfilling this condition.

Suppose that the Gorini-Kossakowski-Sudarshan-Lindblad generator
\beqn{
L(X)
~:=~
ee^*Xee^*-\frac{ee^*X+Xee^*}{2}+i\SB{H,X}
~=~
e\AB{e,Xe}e^*-\frac{e(X^*e)^*+(Xe)e^*}{2}+i\SB{H,X}
}\eeqn
leaves a commutative \nbd{*}subalgebra $\cC$ of $\sB(G)$ invariant. We will show in that case $\dim\cC\le2$. In other words, for every triple $G,e,H$ where $d\ge3$, $L$ does not leave any masa invariant.

Since $L$ leaves also $\cC''$ invariant, we may assume that $\cC$ is a von Neumann algebra, hence, generated by its projections.

Suppose $q_1$ and $q_2$ are two mutually orthogonal projections in $\sB(G)$ such that $q_j$ commutes with $L(q_j)$ for $j=1,2$. We get
\bmu{\label{q12}
0
~=~
q_2L(q_1)q_1
~=~
q_2\,\Bfam{\,e\AB{e,q_1e}e^*-\frac{e(q_1e)^*+(q_1e)e^*}{2}+i\SB{H,q_1}\,}\,q_1
\\
~=~
(q_2e)(q_1e)^*\bfam{\AB{e,q_1e}-\textstyle\frac{1}{2}}
+
iq_2Hq_1,
}\emu

Suppose $q$ is a projection commuting with $L(q)$ such that $qe=0$. Put $q_1:=q$ and $q_2:=\U-q$. Then \eqref{q12} reads $Hq=qHq$. Together with the adjoint equation we get $qH=Hq$. Then, $\CB{H}'\cap\CB{ee^*}'\ni q=\lambda\U$ for some $\lambda\in\C$. By $qe=0$, it follows $q\ne\U$, so, $q=0$. In other words, for every nonzero projection $q$ commuting with $L(q)$ we have $qe\ne0$.

Suppose $q_1$ and $q_2$ are two mutually orthogonal nonzero projections in $\sB(G)$ such that $q_j$ commutes with $L(q_j)$ for $j=1,2$. Exchanging $1$ and $2$ in \eqref{q12} and taking the adjoint, we find
\beqn{
0
~=~
(q_2e)(q_1e)^*\bfam{\AB{e,q_2e}-\textstyle\frac{1}{2}}
-
iq_2Hq_1.
}\eeqn
Summing the two, we get
\beqn{
0
~=~
(q_2e)(q_1e)^*\bfam{\AB{e,(q_1+q_2)e}-1}.
}\eeqn
Since $q_1e\ne0\ne q_2e$, we must have $\AB{e,(q_1+q_2)e}=1$. Since the projection $q:=\U-q_1-q_2$ commutes with $L(q)$ and fulfills $qe=0$, it follows $q=0$. We conclude that every commutative (unital) von Neumann subalgebra that is invariant for $L$, is at most \nbd{2}dimensional. Consequently, if $\dim G\ge3$, then there is no masa invariant for $L$.
\eex

\bex\label{M3noninvex'}
We now seek a unital CP-map $T$ without any invariant masa. For this example we assume that $G=\C^d$ is finite-dimensional. The idea is similar. Start with $X\mapsto ee^*Xee^*$ for some unit vector $e\in G$ and add some \it{Hamiltonian perturbation}. Just that the perturbation should now be in \it{integrated form}, that is, $X\mapsto UXU^*$ for some unitary $U\in M_d$. Also something must be done to normalize $T$ suitably. We take another unit vector $f$ and define
\beqn{
T(X)
~:=~
e\AB{e,Xe}e^*+\AB{f,Xf}(\U-ee^*)+UXU^*
~=~
e\bfam{\AB{e,Xe}-\AB{f,Xf}}e^*+\AB{f,Xf}+UXU^*,
}\eeqn
so that $\frac{T}{2}$ is a unital CP-map. The conditions we pose on $e$, $f$, and $U$, are as follows.
\begin{enumerate}
\item\label{ex1}
For every eigenvector $u$ of $U$ we have $\abs{\AB{e,u}}^2\ne\abs{\AB{f,u}}^2$ and $\AB{e,u}\ne0$.

\item\label{ex2}
$\AB{v,Uv}\ne0$ for all $0\ne v\in G$. (For instance, $\Re U>0$.)
\end{enumerate}
$e$, $f$, and $U$ fulfilling these conditions, obviously, exist in all finite dimensions (and also, when $G$ is infinite-dimensional and separable).

We choose a unit vector $v\in G$ and compute
\beqn{
T(vv^*)
~=~
e\bfam{\abs{\AB{e,v}}^2-\abs{\AB{f,v}}^2}e^*+\abs{\AB{f,v}}^2+(Uv)(Uv)^*.
}\eeqn
If further $vv^*$ commutes with $T(vv^*)$, then
\beqn{
v\AB{v,T(vv^*)v}
~=~
T(vv^*)v
~=~
e\bfam{\abs{\AB{e,v}}^2-\abs{\AB{f,v}}^2}\AB{e,v}+v\abs{\AB{f,v}}^2+Uv\AB{Uv,v}.
}\eeqn
Suppose $\bfam{\abs{\AB{e,v}}^2-\abs{\AB{f,v}}^2}\AB{e,v}=0$, so that $\abs{\AB{e,v}}^2=\abs{\AB{f,v}}^2$ or $\AB{e,v}=0$. Since $\AB{Uv,v}\ne0$, it follows that $v$ is an eigenvector of $U$, contradicting $\abs{\AB{e,v}}^2\ne\abs{\AB{f,v}}^2$ and $\AB{e,v}\ne0$ for every eigenvector of $U$. Therefore, $\abs{\AB{e,v}}^2\ne\abs{\AB{f,v}}^2$ and $\AB{e,v}\ne0$ for every unit vector $v$ such that $vv^*$ commutes with $T(vv^*)$.

Suppose $V=\bCB{v_1,\ldots,v_d}$ is an ONB for $G$ such that the masa $\C v_1v_1^*+\ldots+\C v_dv_d^*$ is invariant for $T$. Fix a $v\in V$. Then for each $i\ne j$ we have
\beqn{
0
~=~
\AB{v_i,T(vv^*)v_j}
~=~
\AB{v_i,e}\bfam{\abs{\AB{e,v}}^2-\abs{\AB{f,v}}^2}\AB{e,v_j}+\AB{v_i,Uv}\AB{Uv,v_j}.
}\eeqn
Since the left summand is nonzero, $\AB{v_i,Uv}\ne0$ for all $i$. We find
\beqn{
\frac{\AB{v_i,e}}{\AB{v_i,Uv}}
~=~
c\frac{\AB{Uv,v_j}}{\AB{e,v_j}}
}\eeqn
for some constant $c\ne0$ and all $i\ne j$. As soon as $d\ge3$, for $i\ne j$ we may choose $k$ such that $i\ne k\ne j$. Then
\beqn{
\frac{\AB{v_i,e}}{\AB{v_i,Uv}}
~=~
c\frac{\AB{Uv,v_k}}{\AB{e,v_k}}
~=~
\frac{\AB{v_j,e}}{\AB{v_j,Uv}}.
}\eeqn
In other words, $\frac{\AB{v_i,e}}{\AB{v_i,Uv}}$ is constant for all $i$, so that $Uv$ is a multiple of $e$. Since $v\in V$ was arbitrary, we find that all $v\in V$ are multiples of $U^*e$. This contradicts unitarity of $U$. Consequently, for $d\ge3$ there is no invariant masa for $T$.
\eex

\brem
Observe that Example \ref{M3noninvex'} gives a unital CP-map without invariant masa, which, by Observation \ref{mapSGob}, gives also rise to a Markov semigroup without invariant masa. But the example is finite-dimensional. Example \ref{M3noninvex} gives a Markov semigroup without invariant masa also when the dimension is countable infinite. But, we do not know whether this means that a single member of the semigroup does not admit masas; see the note following Observation \ref{mapSGob}. So, an example of a unital CP-map without masa in infinite dimension is still missing.
\erem

\section*{Appendix}
\renewcommand{\theemp}{A.\arabic{emp}}
\setcounter{emp}{0}

For the proofs in this appendix we do not make any attempt to be self-contained. Instead, we assume that the reader is familiar with the notions as introduced in Fagnola and Skeide \cite[Section 2]{FaSk07} for the proofs of Theorems \ref{Rebcor} and \ref{Lincor}, plus the necessary notions from Barreto, Bhat, Liebscher and Skeide \cite{BBLS04,LiSk01} about morphisms of time ordered product systems for the proof of Theorem \ref{Linchthm}. Theorems \ref{Rebcor} and \ref{Lincor} are versions specialized to $\sB(G)$ of the results \cite[Theorem 3.1 and 4.2]{FaSk07} for general von Neumann algebras $\cB\subset\sB(G)$. As the intuition of the proof of necessity in the latter results cannot be grasped without a good portion of experience with Hilbert modules, it appears useless to produce a proof for $\sB(G)$, independent of \cite{FaSk07}, that would not even approximately reveal why it works and where it comes from. Just recall that \it{correspondence} it the fashionable term for Hilbert bimodule.

The following result from \cite{FaSk07} about invariance of a maximal commutative subalgebra under CP-maps for general von Neumann algebras is just \cite[Theorem 3.1]{FaSk07} supplemented by the statement in \cite[Observation 3.3]{FaSk07}.

\bitemp[{\cite[Theorem 3.1]{FaSk07}.~}]\label{TcharCCPthm}
Let $\cB\subset\sB(G)$ be a von Neumann algebra on the Hilbert space $G$ and let $T$ be a normal CP-map $T$ on $\cB$. Suppose $E$ is a von Neumann correspondence over $\cB$ and $\xi\in E$ one of its elements such that $T(b)=\AB{\xi,b\xi}$. Furthermore, let $\cC\ni\id_G$ be a maximal commutative von Neumann subalgebra of $\cB$.

Then $T$ leaves $\cC$ invariant, if and only if there exists a \nbd{*}map $\alpha\colon\cC\rightarrow\sB^a(E)$ fulfilling the following properties:
\begin{enumerate}
\item\label{Ta}
The range of $\alpha$ commutes with the left action of elements of $\cC$ on $E$, that is, for all $c_1,c_2\in\cC$ and $x\in E$ we have
\beqn{
c_1\alpha(c_2)x
~=~
\alpha(c_2)c_1x.
}\eeqn

\item\label{Tb}
For all $c\in\cC$ we have
\beqn{
\alpha(c)\xi
~=~
c\xi-\xi c.
}\eeqn
\end{enumerate}
\eitemp

\brem
Every normal CP-map  on a von Neumann algebra can be obtained in that way. For people who like modules: Do the GNS-construction to obtain a correspondence $E_0$ over $\cB$ with a cyclic vector $\xi\in E_0$ having the correct matrix elements; see \cite[Section 2.1]{FaSk07}. Then close $E_0$ suitably to obtain a von Neumann correspondence $E$ following the procedure from Skeide \cite{Ske00b} as explained in \cite[Section 2.3]{FaSk07}. For people who like the classical approach: Do the Stinespring construction \cite{Sti55} to obtain a Hilbert space $H$ with a nondegenerate normal representation $\pi$ of $\cB$ and a map $\xi\in\sB(G,H)$ such that $T(b)=\xi^*\pi(b)\xi$; see \cite[Section 2.2]{FaSk07}. The GNS-module is, then, the strong closure in $\sB(B,H)$ of $\ls\pi(\cB)\xi\cB$; see \cite[Section 2.3]{FaSk07}.
\erem

As we need the same argument in the proof of Theorem \ref{Lincor}, we repeat from \cite{FaSk07} the reduction of Theorem \ref{Rebcor} to Theorem \ref{TcharCCPthm}.

\proof[Proof of Theorem \ref{Rebcor}.~]
If $\cB=\sB(G)$, then $E=\sB(G,G\otimes\eH)$; see \cite[Section 2.4]{FaSk07}. Let $\bfam{e_i}_{i\in I}$ denote an ONB of $\eH$. The family $\bfam{\id_G\otimes e_i}_{i\in I}$ (where $\id_G\otimes e_i$ denotes the mapping $g\mapsto g\otimes e_i$) is, then, an ONB of $E$ in the obvious sense. (See \cite{Ske00b} for quasi ONBs.) Denote by $L_i:=\AB{\id_G\otimes e_i,\xi}$ the coefficients of $\xi$ with respect to this ONB. Then
\beqn{
T(b)
~=~
\sum_{i\in I}L_i^*bL_i
}\eeqn
is a Kraus decomposition of the CP-map $T$ on $\sB(G)$; see \cite[Section 2.4]{FaSk07}. Moreover, every Kraus decomposition can be obtained in that way. (Simply take $\eH:=\C^{\#I}$ with the canonical ONB and define $\xi:=\sum_{i\in I}L_i\otimes e_i$.) The correspondence between maps $\alpha\colon\cC\rightarrow\sB^a(E)$ fulfilling the hypothesis of Theorem \ref{TcharCCPthm} and coefficients $c_{ij}(c)\in\cC$ fulfilling the hypothesis of Theorem \ref{Rebcor} is, then, given by
\beqn{
c_{ij}(c)
~:=~
\bAB{(\id_G\otimes e_i),\alpha(c)(\id_G\otimes e_j)}.
}\eeqn
(Note that the conditions of Theorem \ref{Rebcor} are, clearly, sufficient.  Therefore an $\alpha$ exists and it is easy to see that $\alpha(c)$ can be chosen to have the expansion coefficients $c_{ij}(c)$.)\qed

\bob\label{lincob}
Suppose $T(X)=\sum_{i=I}L_i^*XL_i=\sum_{j\in J}K_j^*XK_j$. Put $\eH:=\C^{\#I}$ and $\eG:=\C^{\#J}$ and denote by $\bfam{e_i}_{i\in I}$ and $\bfam{f_j}_{j\in J}$, respectively, their canonical ONBs. Then $T=\AB{\xi,\bullet\xi}=\AB{\zeta,\bullet\zeta}$ for the elements $\xi:=\sum_{i\in I}L_i\otimes e_i$ and $\zeta:=\sum_{j\in J}K_j\otimes f_j$ of the von Neumann \nbd{\sB(G)}correspondences $E:=\sB(G,G\otimes\eH)$ and $F:=\sB(G,G\otimes\eG)$, respectively. It follows that $v\xi=\zeta$ defines a unique partial isometry $v\in\sB^{a,bil}(E,F)$ that vanishes on $(\sB(G)\xi\sB(G))^\perp$, whose adjoint sends $\zeta$ to $v^*\zeta=\xi$ and vanishes on $(\sB(G)\zeta\sB(G))^\perp$. The superscript $^{bil}$ refers to that the operators are bilinear, that is, they commute with the action of $\sB(G)$. It follows that $v$ must have the form $v=\id_G\otimes V\in\sB(G\otimes\eH,G\otimes\eG)=\sB^a(E,F)$ for some partial isometry $V\in\sB(\eH,\eG)$. If $v_{ji}$ are the matrix elements of $V$ with respect to the canonical ONBs, we find that $K_j=\sum_{i\in I}v_{ji}L_i$ and $L_i=\sum_{j\in J}v_{ji}^*K_j$. We see that the (strongly closed) linear hull is invariant under the choice of the Kraus decomposition. Moreover, $v$ is injective, if and only if $\xi$ generates $E$, and $v$ is surjective, if and only if $\zeta$ generates $F$. If $v$ is bijective, so that it is unitary, then the dimensions of $\eH$ and $\eG$ must coincide, and no Kraus decomposition can have fewer summands than that \hl{minimal} dimension.
\eob

We now quote the criterion for the (bounded) generators of normal CP-semigroups.

\bitemp[{\cite[Theorem 4.2]{FaSk07}.~}]\label{LcharCCPthm}
Let $\cB\subset\sB(G)$ be a von Neumann algebra on the Hilbert space $G$ and let $L$ be a (bounded) normal CCP-map on $\cB$. Suppose $E$ is a von Neumann correspondence over $\cB$ and $d\colon\cB\rightarrow E$ a bounded derivation such that
\beqn{
\AB{d(b),d(b')}
~=~
L(b^*b')-b^*L(b')-L(b^*)b'+b^*L(\U)b'.
}\eeqn
Furthermore, let $\cC\ni\id_G$ be a maximal commutative von Neumann subalgebra of $\cB$.

Then $L$ leaves $\cC$ invariant, if and only if there exist an element $\zeta\in E$ that reproduces $d\upharpoonright\cC$ as
\beqn{
d(c)
~=~
c\zeta-\zeta c,
}\eeqn
a \nbd{*}map  $\alpha\colon\cC\rightarrow\sB^a(E)$ and a self-adjoint element $\gamma\in\cC$ such that the following conditions are satisfied:
\begin{enumerate}
\item\label{L2a}
The range of $\alpha$ commutes with the left action of elements of $\cC$ on $E$, that is, for all $c_1,c_2\in\cC$ and $x\in E$ we have
\beqn{
c_1\alpha(c_2)x
~=~
\alpha(c_2)c_1x.
}\eeqn

\item\label{L2b}
For all $c\in\cC$ we have
\beqn{
\alpha(c)\zeta
~=~
c\zeta-\zeta c.
}\eeqn

\item\label{L1}
For all $c\in\cC$ we have
\beqn{
L(c)-\AB{\zeta,c\zeta}
~=~
\gamma c.
}\eeqn
\end{enumerate}
\eitemp

\proof[Proof of Theorem \ref{Lincor}.~]
Let $L$ be given in the Gorini-Kossakowski-Sudarshan-Lindblad form as stated in Theorem \ref{Lincor} and fix the Hilbert space $\eH:=\C^{\#I}$ with its canonical basis $\bfam{e_i}_{i\in I}$. We observe that if $L$ fulfills the three conditions in Theorem \ref{Lincor}, then by Theorem \ref{Rebcor} applied to the generator with coefficients $K_i:=L_i-c_i$ we see that $L$ leaves $\cC$ invariant.

Suppose now, conversely, that $L$ leaves $\cC$ invariant. By \cite[Sections 2.6 --2.8]{FaSk07} the vector $\xi:=\sum_{i\in I}L_i\otimes e_i$ in the von Neumann \nbd{\sB(G)}correspondence $E:=\sB(G,G\otimes\eH)$ generates a derivation $d(b):=b\xi-\xi b$ that has the required inner products. By Theorem \ref{LcharCCPthm}, there exists a vector $\zeta=\sum_{i\in I}K_i\otimes e_i$ in $E$ such that $d(c)=c\zeta-\zeta c$, that is,
\beqn{
c(\xi-\zeta)
~=~
(\xi-\zeta)c
}\eeqn
for all $c\in\cC$. Therefore, the coefficients $c_i$ of $\xi-\zeta=\sum_{i\in I}(L_i-K_i)\otimes e_i=\sum_{i\in I}c_i\otimes e_i$ must be elements of $\cC$. The rest follows by applying appropriately the other properties that are required in Theorem \ref{Lincor}.\qed

\lf
\proof[Proof of Theorem \ref{Linchthm}.~]
Recall that by \cite{LiSk01} the (continuous) units of the time ordered product system over a von Neumann \nbd{\cB}correspondence $E$ are parameterized as $\xi^\odot(\beta,\xi)=\bfam{\xi_t(\beta,\xi)}_{t\in\R_+}$ where $\beta\in\cB,\xi\in E$. The family of mappings $b\mapsto\AB{\xi_t(\beta,\xi),b\xi_t(\beta,\xi)}$ $(t\in\R_+)$ form a uniformly continuous CP-semigroup with generator $L(b)=\AB{\xi,b\xi}+b\beta+\beta^*b$. In the case of a Gorini-Kossakowski-Sudarshan-Lindblad generator on $\sB(G)$ we have, as in the preceding proofs, $E=\sB(G,G\otimes\eH)$ and $\xi=\sum_{i\in I}L_i\otimes e_i$. The Gorini-Kossakowski-Sudarshan-Lindblad form is minimal, if and only if the single unit $\xi^\odot(\beta,\xi)$ generates the whole time ordered product system. Suppose $F=\sB(G,G\otimes\eG)$ ($\eG=\C^{\#J}$) is another von Neumann \nbd{\sB(G)}correspondence with elements $\alpha\in B,\zeta\in F$ such that $L(b)=\AB{\zeta,b\zeta}+b\alpha+\alpha^*b$. Sending $\xi_t(\beta,\xi)$ to $\xi_t(\alpha,\zeta)$ defines, then, an isometric morphism from the time ordered product system over $E$ into that over $F$. By \cite[Theorem 5.2.1]{BBLS04} morphisms are parameterized by matrices
$\sMatrix{\gamma&\eta^*\\\eta'&a}\in\sB^{a,bil}(\sB(G)\oplus E,\sB(G)\oplus F)$ such that the parameters of the units transform as
\beqn{
(\beta,\xi)
~\longmapsto~
\bfam{\beta+\gamma+\AB{\eta,\xi}~,\,\eta'+a\xi}.
}\eeqn
By \cite[Corollary 5.2.4]{BBLS04} such a morphism is isometric, if and only if $a$ is isometric, $\eta'$ is arbitrary, $\eta=-a^*\eta'$, and $\gamma=ih-\frac{\AB{\eta',\eta'}}{2}$. Interpreting all this properly in terms of the concrete \nbd{\sB(G)}correspondences  and their elements $\xi,\zeta$, and taking also into account that $\sB^{a,bil}(E,F)=\sB^{a,bil}(\sB(G,G\otimes\eH),\sB(G,G\otimes\eG))=\sB(\eH,\eG)$ (because all elements must commute with $\sB(G)$, the coefficients can just be scalar multiples of $\U$),  gives the statement of Theorem \ref{Linchthm}.\qed

\setlength{\baselineskip}{2.5ex}

\newcommand{\Swap}[2]{#2#1}\newcommand{\Sort}[1]{}
\providecommand{\bysame}{\leavevmode\hbox to3em{\hrulefill}\thinspace}
\providecommand{\MR}{\relax\ifhmode\unskip\space\fi MR }
\providecommand{\MRhref}[2]{%
  \href{http://www.ams.org/mathscinet-getitem?mr=#1}{#2}
}
\providecommand{\href}[2]{#2}

\noindent
B.V.\ Rajarama Bhat: \it{Statistics and Mathematics Unit, Indian Statistical Institute Bangalore, R.\ V.\ College Post, Bangalore 560059, India},
E-mail: \tt{bhat@isibang.ac.in},\\
Homepage: \tt{http://www.isibang.ac.in/Smubang/BHAT/}

\lf\noindent
Franco Fagnola: \it{Dipartimento di Matematica, Politecnico di Milano, Piazza Leonardo da Vinci 32, 20133 Milano, Italy},
E-mail: \tt{franco.fagnola@polimi.it},\\
Homepage: \tt{http://www.mate.polimi.it/$\mtt{\tilde{~}}$qp}

\lf\noindent
Michael Skeide: \it{Dipartimento S.E.G.e S., Università degli Studi del Molise, Via de Sanctis, 86100 Campobasso, Italy},
E-mail: \tt{skeide@math.tu-cottbus.de},\\
Homepage: \tt{http://www.math.tu-cottbus.de/INSTITUT/lswas/\_skeide.html}

\end{document}